\documentclass[11pt]{amsart}
\usepackage{amsmath,amssymb,amsthm}
\usepackage{enumitem}
\usepackage{xcolor}
\usepackage{longtable}
\usepackage{tabularx}
\usepackage{booktabs}
\usepackage{graphicx} % Required for inserting images

\def\E{{\mathbb E}}

\def\R{{\mathbb R}}
\def\P{{\mathbb P}}
\def\var{\varepsilon}

\newtheorem{thm}{Theorem}[section]
\newtheorem{cor}[thm]{Corollary}
\newtheorem{lem}[thm]{Lemma}

\theoremstyle{definition}

\theoremstyle{remark}
\newtheorem{rem}[thm]{Remark}
\newtheorem{ass}[thm]{\bfseries Assumption}

\newcommand{\I}{\mathcal I}

\newcommand{\F}{{\mathcal{F}}}
\newcommand{\rmd}{{\rm d}}

\newcommand{\rmosc}{{\rm osc}}

\begin{document}
\title[Gradient Flow of aggregation-diffusion equations]{Energy Dissipation Along Trajectories and Wasserstein Gradient Flow of Aggregation-Diffusion Equations}

\author{Honglei Chen}
\address{H. Chen: School of Mathematical Sciences, Dalian University of Technology, Dalian
116024, P. R. China}
\email{hcl9712285@gmail.com}

\author{Xuewei Wang}
\address{X. Wang (Corresponding author): School of Mathematical Sciences, Dalian University of Technology, Dalian
116024, P. R. China}
\address{SISSA, Trieste 34136, Italy}
\email{Wangxueweii@163.com, xuwang@sissa.it}

\date{\today}



\begin{abstract}
In this work, we established the Wasserstein gradient flow structure of nonlinear Fokker-Planck equations on bounded domains involving both nonlinear diffusion and nonlocal interaction effects. Our approach is based on stochastic analysis. More precisely, we investigated the associated stochastic differential equation (SDE) with normal reflection at the boundary from a trajectorial perspective and derived a semimartingale decomposition of the corresponding energy process. This decomposition leads to a trajectorial energy dissipation identity, which, together with the characterization of the metric derivative, allows us to establish the Wasserstein gradient flow structure of the aggregation-diffusion equation.
\end{abstract}

\maketitle

\section{Introduction}

Nonlocal interactions frequently occur in nature, for instance, to describe the interactions between individuals in a biological population or to describe the gravitational forces between microscopic particles. This setting includes several important examples that have been extensively studied in the literature \cite{ABFS, BDZ}.
In this paper, we consider the following aggregation-diffusion equation 
%nonlinear Fokker-Planck equation with no-flux boundary conditions
\begin{equation}\label{eq:1.1}
 \begin{cases}
 \partial_t p_t(x)
 = \Delta f(p_t(x)) +\nabla\cdot \big(p_t(x)\nabla\left(V+W\ast p_t(x)\right)
 \big)
~ &\text{in }\Omega,\\
 \big[\nabla f(p_t(x))+p_t(x)\nabla\left(V+W\ast p_t(x)\right)\big]\cdot n=0
 &\text{on }\partial\Omega.
 \end{cases}
\end{equation}
where $\Omega\in \R^d$ is an open connected bounded domain and $n$ is the outward normal to the boundary $\partial\Omega$. 
The initial condition is given by a probability density function $p_0$. The smooth function $f:[0, +\infty) \rightarrow \R$ reflects the influence of nonlinearity, and 
smooth functions $V, W:\Omega \rightarrow [0,+\infty)$ play the roles of confinement and interaction potentials, respectively. The convolution is defined by
$W*p_t(x):=\int_{\Omega}W(x-y)p_t(y)\rmd y$.

The Wasserstein gradient flow formulation of the Fokker-Planck equation can be traced back to \cite{JKO}, where the authors established this structure using a time-discretization scheme. Similar approaches have subsequently been applied to various Fokker-Planck equations. In \cite{P} and \cite{BCH}, the authors employed a BDF2-approach and a fully discrete, implicit-in-time finite-volume scheme, respectively, to establish the Wasserstein gradient flow structure of the nonlinear Fokker-Planck equation \eqref{eq:1.1}. Furthermore, \cite{CT} established the slow diffusion limit in the case of $f(r)=r^m$.

However, these results were all established from a macroscopic perspective. In 2022, Karatzas et al. \cite{KST} developed a trajectorial approach based on the microscopic stochastic dynamics associated with the Fokker-Planck equation. Similar methods have also been applied to various other models (see, for instance, \cite{TC, KC}). In the present work, we consider the following McKean-Vlasov SDE
\begin{equation}
\label{eq:reflected-SDE}
\begin{cases}
\begin{aligned}
X_t={}&X_0-\int_0^t
\left[
\nabla V(X_s)+\nabla W*p_s(X_s)
\right]\rmd s\\
&+\int_0^t
\sqrt{\frac{2f(p_s(X_s))}{p_s(X_s)}}\rmd B_s
-\int_0^t n(X_s)\rmd L_s,
\end{aligned}\\[2ex]
X_0=\xi,\qquad X_t\in\overline{\Omega},\\[2ex]
\begin{aligned}
L_t={}&\int_0^t
\mathbf{1}_{\{X_s\in\partial\Omega\}}\rmd L_s,
\qquad L_0=0
\end{aligned}
\end{cases}
\end{equation}
on $\overline{\Omega}$, where $(B_t)_{t\geq 0}$ denotes a standard Brownian motion on a probability space $(\Omega, \mathcal{G}, \P)$ with normal filtration $(\mathcal{G}_t)_{t\geq 0}$ and $\xi$ has density $p_0$. We use the notation $p_t(X_t):=\frac{\rmd P_t}{\rmd x}(X_t)$, $t \geq 0$ for $P_t:={\rm Law}(X_t)$.

We define the free energy functional by
\begin{align}\label{def-F}
\F(p_t):=\int_{\Omega} \Psi(p_t(x))\rmd x+
\int_{\Omega} V(x)p_t(x)\rmd x+
\frac12 \int_{\Omega} W*p_t(x)p_t(x)\rmd x,
\end{align}
where 
\begin{align}\label{h-Phi}
    \psi(u):=\int_1^u \frac{f'(s)}{s}\rmd s,~
    \Psi(u):=\int_0^u \psi(s)\rmd s.
\end{align}
It consists of the internal energy, the potential energy, and the interaction energy.
The corresponding energy dissipation functional is defined by
\begin{equation}\label{def-I}
\I(p_t):=
\int_\Omega p_t(x)
\left|
\nabla \psi(p_t(x))+\nabla V(x)
+\nabla W*p_t(x)
\right|^2\rmd x.
\end{equation}
Proposition 4.10 in \cite{CT} establishes the following energy dissipation identity
\begin{equation}\label{energy-dissipation}
\F(p_t)-\F(p_{t_0})
=-\int_{t_0}^t\mathcal I(p_s)\rmd s,
~ 0\le t_0\le t\le T
\end{equation}
for any probability density function $p_t$ on $\overline{\Omega}$ such that the integral is finite. 
In particular,
\[
\frac{\rmd}{\rmd t}\F(p_t)
=-\mathcal I(p_t).
\]
Thus, the energy functional is nonincreasing along the evolution of the probability densities. 
%The identity \eqref{energy-dissipation} describes energy dissipation at the level of probability distributions. 

From a microscopic perspective, we analyze the SDE \eqref{eq:reflected-SDE} by applying It\^o's formula to the energy process 
\begin{equation}\label{def-v}
\begin{aligned}
\theta_t(x):
={}&\psi(p_t(x))+V(x)+W*p_t(x)\\
&-\int_\Omega f(p_t(y))\rmd y
-\frac12\int_\Omega
W*p_t(y)p_t(y)\rmd y.
\end{aligned}
\end{equation}
which yields
\[
\theta_t(X_t)-\theta_0(X_0)=M_t+F_t.
\]
Furthermore, we obtain the trajectorial version of the energy dissipation identity
\[
\lim_{t\downarrow t_0}
\frac{
\E[\theta_t(X_t)\mid \F_{t_0}]
-\theta_{t_0}(X_{t_0})
}{t-t_0}
=D_{t_0}(X_{t_0}).
\]
The explicit expressions for $M_t$, $F_t$ and $D_t$ are given in Theorem \ref{thm:diss} and Corollary \ref{cor:diss}.
For further details, refer to Subsection \ref{F}.

A comprehensive theory of the dynamical formulation of the Wasserstein metric
\begin{equation}\label{quadratic}
\begin{aligned}
W_2^2(\bar{\mu}_0, \bar{\mu}_1):= \inf \Big\{  &\int_0^1 \!\! \int_{\Omega} |v_t|^2 \rmd \mu_t \rmd t \Big|  \partial_t \mu_t+ {\rm div}\big(v_t \mu_t \big)=0 \text{ holds}\\ & \text{ in the distributional sence, }
\mu_i=\bar{\mu}_i, i=0,1 \Big\}
\end{aligned}
\end{equation}
and the associated gradient flow framework was developed in \cite{AGS}. This framework allows us to characterize the metric derivative of the curve $t \mapsto P_t$. More precisely, we have the Wasserstein metric derivative
\begin{align*}
\lim_{t \downarrow t_0} \dfrac{W_2(P_t, P_{t_0})}{t-t_0}  = (\I(p_{t_0}))^{\frac{1}{2}} ~ \text{ for } 0< t_0 <T.
\end{align*}

By analyzing the perturbed system and combining the results with the metric derivative of the solution curve $t \mapsto P_t$, we show that the solution $p$ evolves in the direction of steepest descent of $\F$ with respect to $W_2$. Consequently, the solution to \eqref{eq:1.1} can be characterized as the Wasserstein gradient flow with the form
\begin{equation*}
\begin{aligned}
\partial_t p_t(x)=&-{\rm grad}_{W_2} \mathcal{F}(p_t(x)) \\
= & {\rm div}\Big( p_t( x) \big(\nabla \psi ( p_t(x) ) +\nabla V(x)+\nabla W*p_t(x) \big) \Big).
\end{aligned}
\end{equation*}

The paper is organized as follows. Section \ref{pre} contains the preliminaries, including notation, assumptions, and regularity.
Section \ref{main results} presents our main results.
In Subsection \ref{F}, we establish the free energy dissipation identity. 
The definition of the Wasserstein metric $W_2$, along with the Wasserstein metric derivative, are given in Subsection \ref{w2}.
In Subsection \ref{GF}, we establish the Wasserstein gradient flow.
In Section \ref{proof}, we provide proofs.

\section{Preliminaries}\label{pre}
In this section, we introduce the notation and standing assumptions used throughout the paper, together with some regular results.

For the reader's convenience, we collect the key symbols used throughout the paper.
If a symbol requires a more precise definition, we provide it when the symbol first appears in the text.

\begin{longtable}{@{}p{0.17\textwidth} p{0.55\textwidth} p{0.15\textwidth}@{}}
\toprule
\textbf{Symbol} & \textbf{Meaning} & \textbf{Reference} \\
\midrule
\endfirsthead

\toprule
\textbf{Symbol} & \textbf{Meaning} & \textbf{Reference} \\
\midrule
\endhead

$T\in (0, +\infty)$ & a terminal time & \\
$P, p$ & measure and its density for the system & \\
%$\mathcal{P}(\Omega)$ & space of probability measures on $\Omega$ & \\
$\mathcal{P}_2(\Omega)$ & space of probability measures on $\Omega$ with finite second moment & \\
$\F$ & free energy functionals &  \eqref{def-F} \\
$\I$ & energy dissipation functionals & \eqref{def-I} \\
$\theta$ & energy processes &  \eqref{def-v}  \\
$D$ & energy dissipation terms of trajectories & \eqref{def-D}  \\
$W_2$ & 2-Wasserstein metric & \eqref{quadratic} \\
$v$ & vector fields & \eqref{vector} \\
${\rm grad}_{W_2} \mathcal{F}$ & Wasserstein gradient flows & \eqref{gradient} \\
$\big|\partial \mathcal{F} \big|_{W_2}$ & Wasserstein metric slopes & Cor. \ref{wslope} \\
\bottomrule
\end{longtable}

Let us define $\psi: (0,+\infty) \rightarrow \R$ and $\Psi: [0,+\infty) \rightarrow \R$ by
\[
    \psi(u)=\int_1^u \frac{f'(s)}{s}\rmd s, ~ \Psi(u)=\int_0^u \psi(s)\rmd s.
\]
Then integration by parts and $f(0)=0$ give
\begin{equation}\label{Phi-identity}
\Psi(u)
=u\psi(u)-\int_0^u s \psi'(s)\rmd s
=u\psi(u)-f(u),~ u>0.
\end{equation}

To ensure the required regularity, we impose the following assumptions.

\begin{ass}[Regularity conditions]\label{reg} \quad
\begin{enumerate}[label=(H\arabic*)]
  \item \label{O} The domain $\Omega$ is an open connected bounded subset of $\R^d$ for some $d\geq 2$, and the boundary $\partial \Omega$ is smooth. 
  \item \label{initial} The initial condition satisfies $p_0\in C^{\infty}(\Omega)$ and the no-flux condition
$$\big[\nabla f(p_0(x)) +p_0\nabla(V(x)+W\ast p_0(x) ) \big]\cdot n=0~ \text{ for all }x \in \partial\Omega.$$
The initial free energy is finite, i.e. $\mathcal{F}(p_0)<+\infty$. Moreover, there exist constants $0< m_0 \leq M_0 <+\infty$ such that 
$$m_0 \leq p_0(x) \leq M_0~ \text{ for all }x \in \Omega. $$
  \item \label{f} The smooth function $f: [0, +\infty) \rightarrow \R$ satisfies 
$$f(0)=f'(0)=0,~ f'(r)>r~ \text{ for all } r>0.$$
  \item \label{potential} The smooth potential functions $V, W: \Omega \rightarrow [0, +\infty)$ are assumed to satisfy that $W$ is symmetric, i.e.,
  \[W(x)=W(-x),~ x\in \R^d.\]
  \item \label{C*} There exist constants $a,b$ such that $0<a <m_0\leq M_0\leq b<+\infty$ and 
  \begin{equation*}\label{H6}
  \begin{aligned}
  \rmosc_{\overline{\Omega}}V+\rmosc_{\overline{\Omega}-\overline{\Omega}}W&+\dfrac{3b|\Omega|\beta_+(b-a)\sqrt{C_P}}{a \beta_-}\|\nabla W\|_{\infty} \\
  & < \beta_- \min\{m_0-a, b-M_0\},
  \end{aligned}
  \end{equation*}
  where $\beta_+:=\max_{[a,b]}\dfrac{f'(r)}{r}$, $\beta_-:=\min_{[a,b]}\dfrac{f'(r)}{r}$,
  $\rmosc_{\overline{\Omega}}V := \sup_{\overline{\Omega}}V - \inf_{\overline{\Omega}} V$,
  $\rmosc_{\overline{\Omega}-\overline{\Omega}}W := \sup_{\overline{\Omega}-\overline{\Omega}}W - \inf_{\overline{\Omega}-\overline{\Omega}} W$
  and $C_P$ denotes the Poincar\'e constant of the domain $\Omega$, i.e.
  $$\Big\|w-\frac{1}{|\Omega|}\int_{\Omega} w \rmd x \Big\|_{L^2(\Omega)}^2 \leq C_P \|\nabla w \|_{L^2(\Omega)}^2 ~\text { for all } w\in H^1(\Omega).$$
  For brevity, we write $C_*:=\dfrac{3b|\Omega|\beta_+(b-a)\sqrt{C_P}}{a \beta_-}$, $\sigma:=\rmosc_{\overline{\Omega}}V+\rmosc_{\overline{\Omega}-\overline{\Omega}}W$ and $k:=\|\nabla W\|_{\infty}$.
\end{enumerate}
\end{ass}

\begin{rem}
(i)By Assumption \ref{f}, $f'(0)=0$ and $f''$ is continuous.
The mean value theorem gives
\[
\left|\frac{f'(s)}{s}\right|
=\left|\frac{f'(s)-f'(0)}{s}\right|
\leq \sup_{r\in[0,1]}|f''(r)|<\infty,
\qquad 0<s\leq1.
\]
Thus the integral defining $\psi$ in \eqref{h-Phi} converges
as $u\downarrow0$, and
\[
\psi(0):=\lim_{u\downarrow0}\psi(u)
=-\int_0^1\frac{f'(s)}{s}\,ds
\]
is finite. Hence $\psi$ extends continuously to $[0,\infty)$,
and $\psi\in L^1_{\mathrm{loc}}([0,\infty))$.\\
(ii) Assumption \ref{C*} ensures that the diffusion effect is stronger than the potential effect, which preventing the formation of vacuum regions. Consequently, the diffusion coefficient in the associated stochastic dynamics remains nondegenerate, a property essential for our subsequent analysis.  \\
(iii) To illustrate the applicability of Assumption \ref{reg}, we provide the following example that satisfies it.
Let $\Omega\in \R^2$ be the open ball of radius $R=\frac{1}{\sqrt{3\pi}}$ centered at the origin, and set
\[
    f(r)=r^3,~
    W(z)=- \frac13\exp{\big(-\frac{|z|^2}{\pi}\big)},~
    V(x)=-\frac{\lambda}{2}|x|^2,
\]
where 
\[
    p_0(x)=3+\frac65\cos(\frac{\pi|x|^2}{R^2}),~
    \lambda=\frac1R\partial_r(W*p_0)(R).
\]
Then 
\[
    \Psi''(r)=\frac{f'(r)}r=3r,~
    \Psi'(r)=\frac32r^2,~
    \frac{f(r)}r=r^2.
\]
Choose
\[
    a=\frac65,~
    m_0=\frac95,~
    M_0=\frac{21}{5},~
    b=\frac{24}{5}.
\]
A direct computation shows that Assumption \ref{reg} is satisfied. In particular, we have 
\[
|\nabla\cdot\bigl[\lambda x-\nabla W*p_t(x)\bigr]|
\le\frac{4}{3\pi}.
\]
\end{rem}

These assumptions ensure the required regularity of the solutions to \eqref{eq:1.1} and \eqref{eq:reflected-SDE}.

\begin{lem}\label{PDEsolution}
Suppose that Assumption \ref{reg} holds. Then the nonlinear Fokker-Planck equation \eqref{eq:1.1} has a unique classical solution 
\begin{align}\label{inte-p}
p_t\in C^{1,2}([0,T] \times \bar{\Omega}),~ \int_{\Omega} p_t(x) \rmd x=1,
\end{align}
\[
a<p_t(x)<b ~\text{ for all } t\in [0,t], x\in \bar{\Omega},
\]
\[
\psi(m_0)-\sigma-C_*k\leq \psi(p_t) \leq \psi(M_0)+\sigma+C_*k. 
\]
Moreover, there exists a unique stationary solution $p_{\infty}$ such that 
\[a<p_{\infty}<b,\]
\begin{equation}\label{state}
\psi(p_{\infty})+V+W*p_{\infty}=C_\infty \text { for some constant } C_\infty
\end{equation}
and 
\begin{equation*}\label{pestimate}
\|p_t-p_{\infty}\|_{L^2(\Omega)} \leq \sqrt{\frac{\beta_+}{\beta_-}}\exp\Big({-\frac{a\beta_-^2 }{2\beta_+C_p}t}\Big)\|p_0-p_{\infty}\|_{L^2(\Omega)}.
\end{equation*}
\end{lem}

\begin{lem}\label{SDEsolution}
Suppose that Assumption \ref{reg} holds. The McKean-Vlasov SDE \eqref{eq:reflected-SDE} has a pathwise unique strong solution $(X,L)$, and the probability density function of $ (X_t)_{t\in [0,T]}$ is the classical solution to the nonlinear Fokker-Planck equation \eqref{eq:1.1}.
\end{lem}

\section{main results}\label{main results}
\subsection{Dissipation of free energy}\label{F}
%Let $p_t$ be the solution of \eqref{eq:1.1} given in Lemma \ref{PDEsolution}, and let $(X,L)$ be the reflected diffusion given in Lemma \ref{SDEsolution}.
In this subsection,
we aim to establish a trajectorial counterpart of the
energy dissipation identity \eqref{energy-dissipation}.

A natural choice of the energy process is
\[
\tilde{\theta}_t(x)
=\frac{\Psi(p_t(x))}{p_t(x)}
+V(x)+\frac12W*p_t(x),
\]
which satisfies
$\E[\tilde{\theta}_t(X_t)]=\F(p_t)$.
However, the no-flux boundary condition yields
\[
\begin{aligned}
\partial_n \theta_t(x)
={}&\left(1-\frac{f(p_t(x))}
{p_t(x)f'(p_t(x))}\right)\partial_n V(x)\\
&+\left(\frac12-\frac{f(p_t(x))}
{p_t(x)f'(p_t(x))}\right)
\partial_n(W*p_t)(x),
\end{aligned}
\]
which does not vanish in general on $\partial\Omega$.
Consequently, applying It\^o's formula to
$\theta_t(X_t)$ produces an additional boundary local
time term. To eliminate this term while preserving
the free energy representation, we choose the 
\emph{energy process} $(\theta_t(X_t))_{t\in [0,T]}$ defined by
\[
\begin{aligned}
\theta_t(x)  
=&\psi(p_t(x))+V(x)+W*p_t(x) \\
&-\int_\Omega f(p_t(y))\rmd y 
-\frac12\int_\Omega
W*p_t(y)p_t(y)\rmd y.
\end{aligned}
\]

We define the \emph{trajectorial energy dissipation
function} $D:[0,T]\times\overline{\Omega}\to\R$ by
\begin{equation}\label{def-D}
\begin{aligned}
D_t(x)
:={}&\partial_t \theta_t(x)
+\frac{f(p_t(x))}{p_t(x)}
\Delta \theta_t(x)\\
&-\big\langle
\nabla \theta_t(x),
\nabla\bigl(V(x)+W*p_t(x)\bigr)
\big\rangle.
\end{aligned}
\end{equation}

\begin{thm}\label{thm:diss}
Suppose that Assumption \ref{reg} holds. Then the energy process $(\theta_t(X_t))_{t\in[0,T]}$ admits the semimartingale decomposition
\begin{align}\label{v-decomposition}
\theta_t(X_t)-\theta_0(X_0)=M_t+F_t ~\text{for } t\in[0,T],
\end{align}
where
\begin{equation}\label{Mt-Ft}
\begin{aligned}
 F_t&:=\int_{0}^{t} D_s(X_s)\rmd s,\\
 M_t&:=\int_{0}^{t} \sqrt{\frac{2f(p_s(X_s))}{p_s(X_s)}}\langle \nabla \theta_s(X_s), \rmd B_s \rangle .
\end{aligned}
\end{equation}
The function $M_t$ is an $L^2$-bounded martingale. Also, we have
\begin{align}\label{EFt}
\E[F_t]=-\int_{0}^{t} \I(p_s)\rmd s>-\infty ~\text{for } t\in[0,T].
\end{align}
\end{thm}

On this basis, we derive the energy dissipation identities at both the macroscopic and trajectorial levels.

\begin{cor}\label{cor:diss}
Suppose that Assumption \ref{reg} holds. Then for every
$0\le t_0\le t\le T$, the energy dissipation identity
\eqref{energy-dissipation} holds. The corresponding
differential version
\begin{equation}\label{eq:diss-rate}
\left.\frac{\rmd}{\rmd t}\F(p_t)\right|_{t=t_0}
=-\I(p_{t_0})
\end{equation}
also holds for $t_0\in[0,T]$. Moreover, for every
$t_0\in[0,T)$, the conditional trajectorial rate of
energy dissipation is given by
\begin{equation}\label{eq:conditional-diss}
\lim_{t\downarrow t_0}
\frac{
\E[\theta_t(X_t)\mid \F_{t_0}]
-\theta_{t_0}(X_{t_0})
}{t-t_0}
=D_{t_0}(X_{t_0}),
\end{equation}
where the limit exists in $L^1(\mathbb P)$.
\end{cor}

\begin{rem}
By \eqref{EFt}, we have
\[
\E[D_{t_0}(X_{t_0})]=-\I(p_{t_0}).
\]
Thus, the identity \eqref{eq:conditional-diss} provides a conditional
trajectorial version of the energy dissipation identity
\eqref{eq:diss-rate}.
\end{rem}

\subsection{Wasserstein metric derivative}\label{w2}

For notational convenience, we set 
\begin{equation}\label{vector}
v_t(x):=-\nabla\big(\psi(p_t(x))+V(x)+W*p_t(x)\big).
\end{equation}
Then \eqref{eq:1.1}  can be rewritten as
\[
\partial_t p_t
+ \nabla \cdot (v_t p_t)=0 ~\text{ in } \Omega.
\]
The Benamou-Brenier theorem provides the following dynamical formulation of the {\em 2-Wasserstein metric}
\begin{align*}
W_2^2(\bar{\mu}_0, \bar{\mu}_1)= \inf \Big\{  &\int_0^1 \!\! \int_{\Omega} |v_t|^2 \rmd \mu_t \rmd t \Big|  \partial_t \mu_t+ {\rm div}\big(v_t \mu_t \big)=0 \text{ holds}\\ & \text{ in the distributional sence, }
\mu_i=\bar{\mu}_i, i=0,1 \Big\}
\end{align*}
for $\bar{\mu}_0, \bar{\mu}_1 \in \mathcal{P}_{ 2}(\Omega)$.
Consequently, the {\em Wasserstein scalar product} at $\mu$ is defined by
\begin{equation}\label{product}
\langle s_1, s_2 \rangle_{\mu} :=\int_{\Omega}  \langle \nabla \xi_1, \nabla \xi_2 \rangle \rmd \mu,
\end{equation}
where $-{\rm div} \big(\mu\nabla \xi_i\big)=s_i$, $i=1,2$.

Definition 8.4.1 of \cite{AGS} provides the following characterization of the tangent bundle
\[
{\rm Tan}_{\mu} \mathcal{P}_2(\Omega):=\overline{\{ \nabla \zeta : \zeta \in C_c^{\infty}(\Omega) \}}^{L^2(\mu)}.
\]
Combining Assumption \ref{O} with the Lemma \ref{PDEsolution}, we conclude that $p$ is an absolutely continuous curve in $(\mathcal{P}_2(\Omega), W_2)$.
Hence, applying Proposition 8.4.3 in \cite{AGS}, we obtain the following result.
\begin{thm}\label{wjder}
Suppose that Assumption \ref{reg} holds and that $v_t\in {\rm Tan}_{P_t} \mathcal{P}_2(\Omega)$. Then the metric derivative of the curve $t \mapsto P_t$ is given by
\begin{equation}\label{mder}
\begin{aligned}
\lim_{t \downarrow t_0} \dfrac{W_2(P_t, P_{t_0})}{t-t_0}  =& \Big(\E[|\nabla\psi(p_{t_0}(X_{t_0}))+\nabla V(X_{t_0})+\nabla W*p_{t_0}(X_{t_0})|^2]\Big)^{\frac{1}{2}} \\
=& (\I(p_{t_0}))^{\frac{1}{2}} ~ \text{ for } 0< t_0 <T.
\end{aligned}
\end{equation}
\end{thm}

\subsection{Wasserstein gradient flow}\label{GF}
We adopt the classical definition given in \cite{AGS}. 
The {\em gradient flow} of $\mathcal{F}$ with respect to $W_2$ at $\bar{\mu}$ is given by 
\begin{equation}\label{gradient}
\Big \langle {\rm grad}_{W_2} \mathcal{F}(\bar{\mu}), \dfrac{\partial \mu_t}{\partial t}\Big|_{t=0} \Big \rangle_{ \bar{\mu}} =\dfrac{\rmd }{\rmd t}\bigg|_{t=0} \mathcal{F}(\mu_t),
\end{equation}
where $\mu: (-t_0, t_0) \rightarrow \mathcal{P}_2(\Omega)$ is an arbitrary smooth curve satisfying $\mu_0=\bar{\mu}$.
The Wasserstein scalar product $\langle \cdot, \cdot \rangle_{\bar{\mu}} $ is given by \eqref{product}.

To show that the solution curve evolves in the direction of steepest descent, we introduce the following Assumption and perturbed system.

\begin{ass}\label{ass-perturbation}
The perturbation potential $\var:\overline{\Omega}\rightarrow\R$ is smooth and $V+\var$ satisfies Assumption \ref{potential}. Moveover, the normal derivative of the perturbation potential vanishes on the boundary, i.e., $\partial_n \var(x)=0$ for $x\in\partial\Omega$.
\end{ass}

Fix $t_0\in[0,T)$. Consider the perturbed equation
\begin{equation}\label{eq:1.2}
 \begin{cases}
 \partial_t p_t^\var
= \Delta f(p_t^\var)
+\nabla\cdot \big(
p_t^\var\nabla\left(V+W\ast p_t^\var+\var\right)
\big) ,
~ &\text{in }\Omega,\\
 \big[
\nabla f(p_t^\var)
+p_t^\var\nabla\left(
V+W\ast p_t^\var+\var
\right)
\big]\cdot n=0
 &\text{on }\partial\Omega, \\
 p_{t_0}^\var(x)= p_{t_0}(x) &\text{in }\overline{\Omega}
 \end{cases}
\end{equation}
and the corresponding SDE with normal reflection
\begin{equation}
\label{eq:reflected-SDE-per}
\begin{cases}
\begin{aligned}
X_t^\var={}&X_{t_0}^\var-\int_{t_0}^t
\left[
\nabla\left(V(X_s^\var)+W*p_s^\var(X_s^\var)+\var(X_s^\var)\right)
\right]\rmd s\\
&+\int_{t_0}^t
\sqrt{\frac{2f(p_s^\var(X_s^\var))}{p_s^\var(X_s^\var)}}\rmd B_s
-\int_{t_0}^t n(X_s^\var)\rmd L_s^\var,
\end{aligned}\\[2ex]
X_{t_0}^\var=X_{t_0},\\[2ex]
\begin{aligned}
L_t^\var={}&\int_{t_0}^t
\mathbf{1}_{\{X_s^\var\in\partial\Omega\}}\rmd L_s^\var,
\qquad L_{t_0}^\var=0.
\end{aligned}
\end{cases}
\end{equation}
where $L^\var$ is continuous and nondecreasing.

\begin{rem}
By Assumption \ref{ass-perturbation} and the boundary
condition in \eqref{eq:1.1}, the initial density $p_{t_0}$
satisfies the boundary compatibility condition in \eqref{eq:1.2}.
Since $V+\var$ is smooth and $p_{t_0}>0$ by
Lemma \ref{PDEsolution}, the local existence argument for
Lemma \ref{PDEsolution} applies to \eqref{eq:1.2}.
Thus, there exists $T_\var>t_0$ such that \eqref{eq:1.2}
has a unique classical solution
\[
p_t^\var\in C^{1,2}([t_0,T_\var]\times\overline{\Omega}), ~
\inf_{[t_0,T_\var]\times\overline{\Omega}}p_t^\var>0, ~
\int_\Omega p_t^\var(x)\rmd x=1.
\]
The argument for Lemma \ref{SDEsolution} then gives a
pathwise unique strong solution $(X^\var,L^\var)$ to
\eqref{eq:reflected-SDE-per}, and $X_t^\var$ has density
$p_t^\var$ for every $t\in[t_0,T_\var]$.
\end{rem}

In parallel with \eqref{def-v}, define the perturbed energy process $(\theta^\var_t(X_t^\var))_{t\in[t_0,T_\var]}$ by
\begin{equation*}\label{def-v-epsilon}
\begin{aligned}
\theta^\var_t(x)
:=&\psi(p_t^\var(x))
+V(x)+W*p_t^\var(x)\\
&-\int_\Omega f(p_t^\var(y))\rmd y
-\frac12\int_\Omega
W*p_t^\var(y)p_t^\var(y)\rmd y.
\end{aligned}
\end{equation*}
We define the corresponding \emph{trajectorial energy dissipation function} $D^\var:[t_0,T_\var]\times\overline{\Omega}\to\R$ by
\begin{equation}\label{def-D-epsilon}
\begin{aligned}
D^\var_t(x)
:={}&\partial_t \theta^\var_t(x)
+\frac{f(p_t^\var(x))}
{p_t^\var(x)}
\Delta \theta^\var_t(x)\\
&-\big\langle
\nabla \theta^\var_t(x),
\nabla\big(
V(x)+W*p_t^\var(x)
+\var(x)\big)
\big\rangle.
\end{aligned}
\end{equation}

\begin{lem}\label{thm:perturbed-diss}
Suppose that Assumptions \ref{reg} and \ref{ass-perturbation} hold. Then the perturbed energy process $(\theta^\var_t(X_t^\var))_{t\in[t_0,T_\var]}$ admits the semimartingale decomposition
\begin{equation}\label{v-epsilon-decomposition}
\theta^\var_t(X_t^\var)
-\theta^\var_{t_0}(X_{t_0}^\var)
=M_t^\var+F_t^\var,
~ t\in[t_0,T_\var],
\end{equation}
where
\begin{equation}\label{Mt-Ft-epsilon}
\begin{aligned}
F_t^\var
&:=\int_{t_0}^t
D^\var_s(X_s^\var)\rmd s,\\
M_t^\var
&:=\int_{t_0}^t
\sqrt{\frac{2f(p_s^\var(X_s^\var))}
{p_s^\var(X_s^\var)}}
\left\langle
\nabla \theta^\var_s(X_s^\var),\rmd B_s
\right\rangle.
\end{aligned}
\end{equation}
The function $M_t^\var$ is an $L^2$-bounded martingale. Moreover, for every $t\in[t_0,T_\var]$,
\begin{equation}\label{EFt-epsilon}
\begin{aligned}
\E[F_t^\var]
={}&-\int_{t_0}^t
\I(p_s^\var)\rmd s\\
&-\int_{t_0}^t
\E\left[
\left\langle
\nabla \theta^\var_s(X_s^\var),
\nabla\var(X_s^\var)
\right\rangle
\right]\rmd s
>-\infty.
\end{aligned}
\end{equation}
\end{lem}

By taking expectations in \eqref{v-epsilon-decomposition}, we obtain the corresponding perturbed energy dissipation identity and its conditional trajectorial version.

\begin{cor}\label{cor:perturbed-diss}
Suppose that Assumptions \ref{reg} and \ref{ass-perturbation} hold. Then, for every $t\in[t_0,T_\var]$,
\begin{equation}\label{eq:perturbed-energy-dissipation}
\begin{aligned}
&\F(p_t^\var)
-\F(p_{t_0}^\var)\\
=&-\int_{t_0}^t
\I(p_s^\var)\rmd s-\int_{t_0}^t
\E\left[
\left\langle
\nabla \theta^\var_s(X_s^\var),
\nabla\var(X_s^\var)
\right\rangle
\right]\rmd s.
\end{aligned}
\end{equation}
The corresponding differential version also holds, i.e.
\begin{equation}\label{eq:perturbed-diss-rate}
\begin{aligned}
\left.
\frac{\rmd^+}{\rmd t}
\F(p_t^\var)
\right|_{t=t_0}
=-\I(p_{t_0})
-\E\left[
\left\langle
\nabla \theta_{t_0}(X_{t_0}),
\nabla\var(X_{t_0})
\right\rangle
\right].
\end{aligned}
\end{equation}
Moreover, the conditional trajectorial rate of energy dissipation is given by
\begin{equation}\label{eq:conditional-perturbed-diss}
\lim_{t\downarrow t_0}
\frac{
\E\left[
\theta^\var_t(X_t^\var)
\mid \F_{t_0}
\right]
-\theta^\var_{t_0}(X_{t_0}^\var)
}{t-t_0}
=D^\var_{t_0}(X_{t_0}),
\end{equation}
where the limit exists in $L^1(\mathbb P)$.
\end{cor}

\begin{rem}
As in the unperturbed case, the expectation of the
trajectorial energy dissipation function satisfies
\[
\E[D^\var_t(X_t^\var)]
=-\I(p_t^\var)
-\E\left[
\left\langle
\nabla \theta^\var_t(X_t^\var),
\nabla\var(X_t^\var)
\right\rangle
\right].
\]
Thus, the expected trajectorial energy dissipation rate
includes an additional contribution from the perturbation.
\end{rem}

Similar with Theorem \ref{wjder}, we have following Lemma. 
\begin{lem}\label{dervar}
Suppose that Assumptions \ref{reg} and \ref{ass-perturbation} hold. The pertubed vector field $v^\var_t:=-\nabla\big(\psi(p^\var_t)+V+W*p^\var_t+\var\big)\in {\rm Tan}_{P^\var_t} \mathcal{P}_2(\Omega)$. Then the metric derivative of the curve $t \mapsto P^\var_t$ is given by
\begin{equation*}\label{mder2}
\begin{aligned}
& \lim_{t \downarrow t_0} \dfrac{W_2(P^\var_t, P^\var_{t_0})}{t-t_0} \\  
=& \Big(\E[|\nabla \psi(p^\var_{t_0}(X^\var_{t_0}))+\nabla V(X^\var_{t_0})+\nabla W*p^\var_{t_0}(X_{t_0})+\nabla\var(X^\var_{t_0})|^2]\Big)^{\frac{1}{2}}
\end{aligned}
\end{equation*}
for $0< t_0 <T$. 
\end{lem}

We now state the Wasserstein gradient flow theorem, together with the corresponding characterization of the Wasserstein metric slope.

\begin{thm}\label{main}
Suppose that Assumption \ref{reg} holds, and that $v_t\in {\rm Tan}_{P_t} \mathcal{P}_2(\Omega)$.
The Wasserstein gradient flow of the free energy functional $\mathcal{F}$ with respect to $W_2$ is
\begin{equation}\label{form2}
\begin{aligned}
\partial_t p_t(x)=&-{\rm grad}_{W_2} \mathcal{F}(p_t(x)) \\
= & {\rm div}\Big( p_t( x) \big(\nabla \psi ( p_t(x) ) +\nabla V(x)+\nabla W*p_t(x) \big) \Big).
\end{aligned}
\end{equation}
\end{thm}

\begin{cor}\label{wslope}
The Wasserstein metric slope of the free energy functional $\mathcal{F}$ along the curves $t \mapsto P_t$ is given by
\begin{equation*}
\begin{aligned}
\big|\partial \mathcal{F} \big|_{W_2}(P_{t_0}):=
&\lim_{t \downarrow t_0} \dfrac{\F(P_t) -\F(P_{t_0})}{W_2\big(P_t, P_{t_0}\big)} \\
=&(\I(p_{t_0}))^{\frac{1}{2}}, ~ \text{ for } 0< t_0 <T.
\end{aligned}
\end{equation*}

\end{cor}

\section{Proof of the results}\label{proof}
\subsection{Proof of Lemma \ref{PDEsolution} and \ref{SDEsolution}} 
\begin{proof}[Proof of Theorem \ref{PDEsolution}]
We first establish the existence of a local-in-time solution and a stationary. We then combine a priori estimate with the maximum principle to derive upper and lower bounds for the solution. Using these bounds, we show that the maximal existence time is infinite. Finally, we establish further regularity and obtain the existence and uniqueness of a global classical solution.

{\bf Step 1.} local-in-time solution and a stationary.

By Theorem 13.1 of \cite{A}, the equation \eqref{eq:1.1} admits a unique maximal weak solution defined on the maximal interval of existence $[0, T_{\max})$. 
Conservation of mass follows by
\[
\frac{\rmd}{\rmd t}\int_\Omega p_t(x)\rmd x
 =\int_{\partial\Omega}\big(\nabla f(p_t)+p_t\nabla(V+W*p_t)\big)\cdot n\rmd S=0.
\]

Since $m_0\leq \bar p:=\frac{1}{|\Omega|}\leq M_0$,
Assumption \ref{C*} gives
\[
\sigma<\beta_-\min\{m_0-a,b-M_0\}
\le \min\{\psi(\bar p)-\psi(a),\psi(b)-\psi(\bar p)\}.
\]
Thus \((\psi|_{[a,b]})^{-1}\) is defined on
\([\psi (\bar p)-\sigma,\psi (\bar p)+\sigma]\subset [\psi (a),\psi (b)]\).
It maps this interval into \((a,b)\).
Denote 
\[
Q:=\left\{\rho\in C(\overline\Omega)\ \middle|\
a\le\rho\le b,~ \int_\Omega\rho\rmd x=1\right\}.
\]
Then \(\bar p\in Q\) and 
\[
\rmosc_{\overline\Omega}(V+W*\rho)
\le \rmosc_{\overline\Omega}V
   +\rmosc_{\overline{\Omega}-\overline{\Omega}}W=\sigma \text{ for } \rho\in Q.
\]
The function \((\psi|_{[a,b]})^{-1}(c-V-W*\rho)\) is well defined for \(\psi (\bar p)+\min_{\overline\Omega}\{V+W*\rho\}\le c\le \psi (\bar p)+\max_{\overline\Omega}\{V+W*\rho\} \),
and its continuity and
strict monotonicity give a unique constant \(C_\rho\) such that the map
\[
\Lambda: \rho \mapsto (\psi|_{[a,b]})^{-1}(C_\rho-V-W*\rho)
\]
satisfies
$\int_\Omega \Lambda\rho\rmd x=1$.
In particular, for every \(x\in\overline\Omega\),
\[
a<(\psi|_{[a,b]})^{-1}(\psi (\bar p)-\sigma)\le (\Lambda\rho)(x)
 \le (\psi|_{[a,b]})^{-1}(\psi (\bar p)+\sigma)<b.
\]
Hence \(\Lambda(Q)\subset Q\).
Applying the Schauder's fixed point theorem (see Theorem 8.1 in \cite{L}) to 
$\Lambda\rho=(\psi|_{[a,b]})^{-1}(C_\rho-V-W*\rho)$ for $\rho\in Q,$ 
we obtain a fixed point $p_{\infty}$, which satisfies 
\[
 a<p_\infty<b,~ \int_\Omega p_\infty\rmd x=1,~
 \psi (p_\infty)+V+W*p_\infty=C_{p_\infty}=C_\infty.
\]
Assumptions \ref{f}, \ref{potential} and the identity $p_\infty=(\psi|_{[a,b]})^{-1}(C_\infty-V-W*p_\infty)$ yield $p_\infty\in C^\infty(\overline\Omega)$.
Differentiating \eqref{state} and using $r\Psi''(r)=f'(r)$ gives
\[
 \nabla f(p_\infty)+p_\infty\nabla(V+W*p_\infty)
 =p_\infty\nabla\big(\psi (p_\infty)+V+W*p_\infty\big)=0 .
\]
Thus $p_\infty$ satisfies the stationary equation and its no-flux
boundary condition.

{\bf Step 2.} A priori estimate.

Before establishing the upper and lower bounds for the solution, we first derive a key estimate. To this end, we consider the local classical solution $p$ on the time interval before it first attains either $a$ or $b$ on $\overline{\Omega}$. The continuity of $p$, together with the strict initial bounds, guarantees the existence of such an interval.

According to the definition of $C_*$ and \ref{C*},
\begin{equation}\label{ineq0}
  \frac{b|\Omega|k\sqrt{C_P}}{\beta_-}
  =\frac{a}{3\beta_+(b-a)}\,C_*k
  \le\frac{a\beta_-\min\{m_0-a,b-M_0\}}
              {3\beta_+(b-a)}
   \le \frac a6,
\end{equation}
where the last inequality follows from
$ 2\min\{m_0-a,b-M_0\}\le (m_0-a)+(b-M_0)\le b-a $.
Denote
\[
  w_t:=p_t-p_\infty,~
  \eta_t:=\psi(p_t)-\psi(p_\infty),~
  \bar\eta:=\frac1{|\Omega|}\int_\Omega\eta\rmd x .
\]
Conservation of mass yields $\int_\Omega w_t\rmd x=1-1=0$. 
Combining $\Psi''\ge \beta_-$ on $[a,b]$, the Cauchy-Schwarz inequality
and the Poincar\'e inequality, we deduce that
\begin{align*}
  \beta_-\|w_t\|_2^2
    &\le\int_\Omega\eta_t w_t\rmd x
     =\int_\Omega(\eta_t-\bar\eta)w_t\rmd x\\
    &\le\|\eta_t-\bar\eta\|_2\|w_t\|_2
     \le\sqrt{C_P}\|\nabla\eta_t\|_2\|w_t\|_2 .
\end{align*}
Consequently,
\begin{equation}\label{w-eta}
  \|w_t\|_2
    \le\frac{\sqrt{C_P}}{\beta_-}\|\nabla\eta_t\|_2.
\end{equation}
Recall that
$\nabla\big(\psi(p_\infty)+V+W*p_\infty\big)=0$.
Hence
\begin{align*}
  &\bigl(\nabla \psi(p_t)+\nabla V+\nabla W*p_t\bigr) \\
   =&\bigl(\nabla \psi(p_t)+\nabla V+\nabla W*p_t\bigr)
     -\bigl(\nabla \psi(p_\infty)+\nabla V+\nabla W*p_\infty\bigr)\\
   =&\nabla\eta_t+\nabla W*w_t .
\end{align*}

Define
\[
  K(t):=\int_\Omega
       \bigl[\Psi(p_t)-\Psi(p_\infty)-\psi(p_\infty)(p_t-p_\infty)\bigr]\rmd x .
\]
Taylor's formula gives
\begin{equation}\label{B}
  \frac{\beta_-}2\|w_t\|_2^2
       \le K(t)\le\frac{\beta_+}2\|w_t\|_2^2 .
\end{equation}
Using $\partial_t p_t=\nabla\cdot\big(p_t\nabla (\psi(p_t)+V+W*p_t)\big)$, we obtain
\begin{align*}
  K'(t)
   =&\int_\Omega \eta_t\,\partial_t p_t\rmd x
     =\int_\Omega\eta_t\,\nabla\cdot\big(p_t\nabla(\psi(p_t)+V+W*p_t)\big)\rmd x\\
   =&-\int_\Omega p_t\nabla\eta_t\cdot\nabla\big(\psi(p_t)+V+W*p_t\big)\rmd x \\
   &+\int_{\partial\Omega}\eta_t\,p_t\nabla\big(\psi(p_t)+V+W*p_t\big)\cdot n\rmd S\\
   =&-\int_\Omega p_t|\nabla\eta_t|^2\rmd x
     -\int_\Omega p_t\nabla\eta_t\cdot\nabla W*w_t\rmd x .
\end{align*}
Using 
\[
\begin{aligned}
  |\nabla W*w_t(x)|
  =& \left|\int_\Omega\nabla W(x-y)w(y)\rmd y\right| \\
  \le& k\int_\Omega|w(y)|\rmd y 
  \le k|\Omega|^{\frac12}\|w_t\|_2 ~ \text{ for } x\in\overline{\Omega}
\end{aligned}
\]
and \eqref{w-eta}, we deduce 
\begin{equation*}\label{conv}
  \|\nabla W*w_t\|_2
   \le k|\Omega|\|w_t\|_2
   \le\frac{k|\Omega|\sqrt{C_P}}{\beta_-}\|\nabla\eta_t\|_2 .
\end{equation*}
It follows that
\[
  \left|\int_\Omega p_t\nabla\eta_t\cdot\nabla W*w_t\rmd x\right|
   \le b\|\nabla\eta_t\|_2\|\nabla W*w_t\|_2
   \le\frac{bk|\Omega|\sqrt{C_P}}{\beta_-}\|\nabla\eta_t\|_2^2
\]
and 
\[
\int_\Omega p_t|\nabla\eta_t|^2\rmd x\ge a\|\nabla\eta_t\|_2^2.
\]
Using \eqref{ineq0}, \eqref{w-eta} and \eqref{B}, we get
\begin{equation}\label{B'}
\begin{aligned}
K'(t)
 \le &-(a-\frac{bk|\Omega|\sqrt{C_P}}{\beta_-}) \|\nabla\eta_t\|_2^2 \\
 \le & -\frac{5a}{6}\|\nabla\eta_t\|_2^2 
 \le  -\frac{5a\beta_-^{2}}{6C_P}\|w_t\|_2^2 
 \le  -\frac{5a\beta_-^{2}} {3\beta_+C_P}K(t).
\end{aligned}
\end{equation}
Therefore, we conclude that
\[
  K(t)\le K(0)\exp\Big(-\frac{5a\beta_-^{2}}{3\beta_+C_P} t\Big)
\]
and 
\[
 \|p_t-p_\infty\|_2^2
    \le\frac{2K(0)}{\beta_-}\exp\Big(-\frac{5a\beta_+^{2}}{3\beta_+C_P} t\Big).
\]
For $K(t)>0$, we obtain
\[
  -\frac{\rmd}{\rmd t}\sqrt{K(t)}
    =-\frac{K'(t)}{2\sqrt{K(t)}}
    \ge\frac{5a}{12}\frac{\|\nabla\eta_t\|_2^2}{\sqrt{K(t)}}
\]
and 
\[
  \sqrt{K(t)}
  \le\frac{\sqrt{\beta_+}}{\sqrt{2}}\|w_t\|_2
  \le\frac{\sqrt{\beta_+C_P}}{\sqrt{2}\beta_-}
       \|\nabla\eta_t\|_2
\]
by \eqref{B'}, \eqref{B} and \eqref{w-eta}.
It follows that
\[
  -\frac{\rmd}{\rmd t}\sqrt{K(t)}
    \ge\frac{5a\sqrt{2}\beta_-}
             {12\sqrt{\beta_+C_P}}\|\nabla\eta_t\|_2 .
\]
Integrating in time yields 
\[
  \int_0^t\|\nabla\eta_s\|_2\rmd s
  \le\frac{12\sqrt{\beta_+C_P}}
           {5a\sqrt{2}\beta_-}\sqrt{K(0)} .
\]
If $K$ vanishes, we integrate up to the first time it reaches zero, which is well-defined by continuity. Since $K$ is nonnegative and nonincreasing, it remains zero thereafter.
Therefore, we have
\begin{equation*}\label{int0t}
  \int_0^t\|\nabla\eta_s\|_2\rmd s
   \le\frac{6\beta_+\sqrt{C_P|\Omega|}} {5a\beta_-}(b-a).
\end{equation*}
by
$K(0)\le\frac{\beta_-}2\|p_0-p_\infty\|_2^2 \le\frac{\beta_+}2|\Omega|(b-a)^2$.

Define $J_t:=p_t\nabla(\psi(p_t)+V+W*p_t)$. The equations $\partial_t p_t=\nabla\cdot J_t$ and $J_t\cdot n=0$
ensure 
\begin{align*}
   \partial_t(W*p_t)(x) 
    =&\int_\Omega W(x-y)\,\nabla_y\cdot J_t(y)\rmd y\\
    =&-\int_\Omega\nabla_y W(x-y)\cdot J_t(y)\rmd y
       +\int_{\partial\Omega}W(x-y)J_t(y)\cdot n(y)\rmd S_y\\
    =&\int_\Omega\nabla W(x-y)\cdot J_t(y)\rmd y.
\end{align*}
From 
\begin{align*}
  \|J_t\|_1
   &\le b\int_\Omega\big(|\nabla\eta_t|+|\nabla W*w_t|\big)\rmd x\\
   &\le b\sqrt{|\Omega|}\big(\|\nabla\eta_t\|_2+\|\nabla W*w_t\|_2\big)\\
   &\le b\sqrt{|\Omega|}
      \left(1+\frac{k|\Omega|\sqrt{C_P}}{\beta_-}\right)
      \|\nabla\eta_t\|_2,
\end{align*}
we can deduce that
$\|\partial_t(W*p_t)\|_\infty\le k\|J_t\|_1$. 
By integrating $\|\partial_t(W*p_t)\|_\infty$ with respect to time from $0$ to $t$, we obtain
\begin{equation*}\label{intW}
\begin{aligned}
  &\int_0^t\|\partial_s(W*p_s)\|_\infty\rmd s \\
    \le& k b\sqrt{|\Omega|}
      \left(1+\frac{k|\Omega|\sqrt{C_P}}{\beta_-}\right)
      \frac{6\beta_+\sqrt{C_P|\Omega|}}
           {5a\beta_-}(b-a)\\
    =&\frac25 C_*k
      \left(1+\frac{k|\Omega|\sqrt{C_P}}{\beta_-}\right)
     \le \frac25 C_*k\left(1+\frac{a}{6b}\right)
     \le C_*k .
\end{aligned}
\end{equation*}

{\bf Step 3.} Upper and lower bounds. 

We continue our analysis on the time interval where $a<p_t<b$. A direct computation yields
\begin{align*}
 & \partial_t(\psi(p_t)+V+W*p_t ) \\
% =&\psi'(p_t )\partial_t p_t +\partial_t(W*p_t )\\
 %=&\frac{f'(p_t )}{p_t} 
  % \bigl(p_t \Delta(\psi(p_t)+V+W*p_t) +\nabla p_t \cdot\nabla(\psi(p_t)+V+W*p_t) \bigr)
  % +\partial_t(W*p_t )\\
 =&f'(p_t )\Delta(\psi(p_t)+V+W*p_t) \\
 &+\nabla \psi(p_t )\cdot\nabla(\psi(p_t)+V+W*p_t) +\partial_t(W*p_t ).
\end{align*}
and the no-flux condition $\partial_n(\psi(p_t)+V+W*p_t)=0$. 

Denote 
\[
 H(t):=\int_0^t\|\partial_s(W*p_s)\|_\infty\rmd s .
\]
Then $0\le H(t)\le C_*k$ and
\[
 -H'(t)\le\partial_t(W*p_t)(x)\le H'(t).
\]
For the function 
\[
 u_t:=(\psi(p_t)+V+W*p_t)
        -\max_{\overline\Omega}\{\psi(p_0)+V+W*p_0\}-H(t),
\]
we have 
\[
 \partial_tu_t-f'(p_t)\Delta u_t-\nabla \psi(p_t)\cdot\nabla u_t
 =\partial_t(W*p_t)-H'(t)\le0,
\]
and 
\[
\partial_nu_t=0,~ u_0\le0.
\]
The maximum principle (see, e.g. Theorem 2.9 in \cite{L}) gives $u_t\le0$,
which means that $(\psi(p_t)+V+W*p_t)(x)\le \max_{\overline\Omega}\{\psi(p_0)+V+W*p_0\}+H(t)$. 
The same argument implies
$(\psi(p_t)+V+W*p_t)(x)\ge \min_{\overline\Omega}\{\psi(p_0)+V+W*p_0\}-H(t).$
We combine 
\[
 \bigl|V(y)+W*p_0(y)-V(x)-W*p_t(x)\bigr|\le \sigma
\]
and the monotonicity of $\psi$ to deduce that
\begin{align*}
 \psi(p_t(x))
 =& (\psi(p_t)+V+W*p_t) (x)-V(x)-W*p_t(x)\\
 \le &\max_{\overline\Omega}\big\{\psi(p_0)+V+W*p_0\big\}+H(t) -V(x)-W*p_t(x)\\
 \le &\psi(M_0) +\max_{\overline\Omega}\big\{V+W*p_0\big\} +H(t)\\
   & -\min_{\overline\Omega}\bigl\{V+W*p_t\bigr\}\\
 \le & \psi(M_0)+\sigma+H(t)
\end{align*}
and 
\begin{align*}
 \psi(p_t(x))
 \ge & \psi(m_0) +\min_{\overline\Omega}\big\{V+W*p_0\big\} -H(t)
    -\max_{\overline\Omega}\big\{V+W*p_t\big\}\\
 \ge & \psi(m_0)-\sigma-H(t).
\end{align*}
Therefore, we have
\[
 \psi(m_0)-\sigma-C_*k
 \le \psi(p_t(x))
 \le \psi(M_0)+\sigma+C_*k.
\]
According to Assumption \ref{C*}, we get
\[
 \psi(m_0)-\psi(a)
 =\int_a^{m_0}\psi'(r)\rmd r
 \ge\beta_-(m_0-a)>\sigma+C_*k
\]
and 
\[
 \psi(b)-\psi(M_0)
 =\int_{M_0}^b \psi'(r)\rmd r
 \ge\beta_-(b-M_0)>\sigma+C_*k.
\]
It follows that
\begin{equation}\label{ub}
 \psi(a)<\psi(m_0)-\sigma-C_*k
 \le \psi(p_t(x))
 \le \psi(M_0)+\sigma+C_*k<\psi(b).
\end{equation}
In fact, the solution cannot reach $a$ or $b$ before $T_{\max}$.
Otherwise, by continuity, the estimate \eqref{ub} would
remain valid at the first contact time, implying that the solution
is still strictly between $a$ and $b$ everywhere in $\overline\Omega$.
This is a contradiction. 
Therefore, we obtain
\[
 a<p_t(x)<b ~\text{ for all }(t,x)\in[0,T_{\max})\times\overline\Omega.
\]

{\bf Step 4.} The maximal existence time is infinite.

Assume that $T_{\max}<+\infty$ and argue by contradiction.
Theorem 6.28 and Theorem 6.44 of \cite{L} imply that for given constants $\alpha\in(0,1)$, $R_\Omega>0$ and $C_0>0$, we have
\[
 \rmosc_{Q_r}p
 \le C_0b\big(1+\|\nabla V\|_\infty+\|\nabla W\|_\infty\big) \left(\frac rR\right)^\alpha,
\]
for $0<r<R\le\min\{R_\Omega,1\}$ and $t-R^2>0$,
where $Q_R:=(\Omega\cap U_R(x))\times(t-R^2,t]$, $U_R$ is a ball of radius $R$.

Fix
$
 \tau=\frac{T_{\max}}2,~
 R=\min\Big\{R_\Omega,1,\frac{\sqrt{\tau}}2\Big\}.
$
Then $t-R^2\ge \frac{3\tau}{4}>0$
and the preceding estimate remains valid. 
Define
\[
 h:=|x-y|+|t-s|^{\frac12},~\text{ for } \tau\le s\le t<T_{\max},~ x,y\in\overline\Omega.
\]
In the case of $0<h<\frac R4$, let $r=2h$.
Then we obtain
\[
 |p_t(x)-p_s(y)|
 \le C_0b(1+\|\nabla V\|_\infty+\|\nabla W\|_\infty) \left(\frac{2h}{R}\right)^\alpha.
\]
If $h\ge \frac R4$, we have
\[
 |p_t(x)-p_s(y)|\le2b
 \le2b\left(\frac4R\right)^\alpha h^\alpha.
\]
We conclude that 
\begin{align*}
 \sup_{[\tau, T_{\max})}\|p_t\|_{C^\alpha(\overline\Omega)}&\le C_\tau\bigl(|x-y|^\alpha+|t-s|^{\frac\alpha2}\bigr)
 <+\infty,\\
  \|p_t-p_s\|_\infty&\le C_\tau|t-s|^{\frac\alpha2},
\end{align*}
where $C_\tau$ is independent of $s,t$.
Since $p_t$ is Cauchy as $t\uparrow T_{\max}$, we define 
\[
p_*:=\lim_{t\uparrow T_{\max}}p_t ~\text{in }C(\overline\Omega)
\]
as the definition of \(p\) at time $T_{\max}$,
which satisfies
\[
 a<p_*<b,~ \int_\Omega p_*\rmd x=1.
\]
It follows that
\[
 p\in C^{\alpha,\frac\alpha2} (\overline\Omega\times[0,T_{\max}]).
\]

We now treat
\(
f'(p_t(x)), 
\nabla V(x)+\nabla W*p_t(x)
\)
as given coefficients and consider the following linear problem
\[
 \begin{cases}
 \partial_tu-\nabla\cdot\big(f'(p)\nabla u+(\nabla V+\nabla W*p)u\big)=0
       &\text{in }\Omega\times(\tau,T_{\max}),\\
 \big(f'(p)\nabla u+(\nabla V+\nabla W*p)u\big)\cdot n=0
       &\text{on }\partial\Omega\times(\tau,T_{\max}),\\
 u_\tau=p_\tau
       &\text{in }\Omega.
 \end{cases}
\]
Then we have
$ u\in C^{1+\alpha,(1+\alpha)/2}
       \bigl(\overline\Omega\times[\tau,T_{\max}]\bigr)$
with the help of Theorem 6.49 in \cite{L}.
Since $u$ and $p$ solve the same linear problem on $[\tau,T], T\in(\tau,T_{\max})$,  with the same initial data, and since $T$ is arbitrary, it follows that
\[
 u=p ~\text{on }\overline\Omega\times[\tau,T_{\max}].
\]
which implies
\[
 p\in C^{1+\alpha,\frac{1+\alpha}2}
       \bigl(\overline\Omega\times[\tau,T_{\max}]\bigr).
\]

We next improve the regularity once more.
Similarly, by treating the coefficients of the nonlinear equation as given, we obtain
\[
 p\in C^{2+\alpha,1+\frac\alpha2} (\overline\Omega\times[0,T_{\max}]).
\]
by Theorem 5.18 of \cite{L}.
We restart the equation from $p_*$ at time $T_{\max}$ and obtain a local classical solution beyond $T_{\max}$.
This contradicts the maximality of $T_{\max}$, and
therefore
\[
 T_{\max}=+\infty.
\]

Global continuation and local classical regularity imply that,
for every finite $T>0$,
\begin{equation*}
\begin{aligned}
p\in C^{2+\alpha,1+\frac\alpha2}
       \big(\overline\Omega\times[0,T]\big),~
 \sup_{[0,T]}\|\nabla p_t\|_\infty<+\infty, \\
 a\le p_t(x)\le b,~
 \int_\Omega p_t(x)\rmd x=1~
 \text{for all }t\ge0.
\end{aligned}
\end{equation*}
Uniqueness follows from Step 1.
Hence the equation admits a unique global classical solution.
\end{proof}

\begin{proof}[Proof of Theorem \ref{SDEsolution}]
Let $p$ be the global classical solution obtained in Lemma \ref{PDEsolution}. 
Appling
\[
    \|\nabla V+\nabla W*p_t\|_{L^\infty(\Omega)}
    \le
    \|D^2V\|_{L^\infty(\Omega)}
    +\|D^2W\|_{L^\infty(\Omega-\Omega)},
    ~ 0\le t\le T
\]
and 
\[
\begin{aligned}
    & \sup_{[0,T]} 
    \Big\|\nabla\sqrt{\frac{2f(p_t)}{p_t}}\Big\|_{L^\infty(\Omega)}\\
    \le&  \sqrt{2}\,
    \Big\|\left(\sqrt{\frac{f(r)}{r}}\right)'\Big\|_{L^\infty([a,b])}
    \sup_{[0,T]}\|\nabla p_t\|_{L^\infty(\Omega)}
    <+\infty
\end{aligned}
\]
to Theorem 1 of \cite{MR}, \eqref{eq:reflected-SDE} has a unique strong solution \((X,L)\) on \([0,T]\). 

We now prove that ${\rm Law}(X_t)=p_t(x)\rmd x$.
Given \(\varphi\in C_c^\infty(\Omega)\), consider the backward problem
\begin{equation}\label{eq:backward}
\begin{cases}
    \partial_t u+\dfrac{f(p)}{p}\Delta u-\nabla (V+W*p)\cdot\nabla u=0
        & \text{in }(0,T)\times\Omega,\\
    \partial_nu=0
        & \text{on }(0,T)\times\partial\Omega,\\
    u_T=\varphi
        & \text{in }x\in\overline\Omega.
\end{cases}
\end{equation}
After reversing time, Theorem 5.18 of \cite{L} guarantees the existence of a classical solution $u$, together with the boundedness of $u$ and all its derivatives required in the subsequent analysis.
Applying It\^o's formula for continuous semimartingales (\cite{KS}, Theorem 3.6) to \(u_t(X_t)\) and using \eqref{eq:backward}, we obtain
\[
    \rmd u_t(X_t)
    =
    \sqrt{\frac{2f(p_t(X_t))}{p_t(X_t)}}
        \nabla u_t(X_t)\cdot \rmd B_t
    -\nabla u_t(X_t)\cdot n(X_t) \rmd L_t.
\]
Applying the same calculation to the solution started from \(x\)
at time \(t\), integrating to \(T\), and taking expectations yields
\[
    u_t(x)=\mathbb{E}_{t,x}\bigl[\varphi(X_T)\bigr].
\]
In particular, for the initial law \(p_0(x)\rmd x\),
\begin{equation}\label{eq:expectation}
    \mathbb{E}\bigl[\varphi(X_T)\bigr]
    =
    \mathbb{E}\bigl[u_0(X_0)\bigr]
    =
    \int_\Omega u_0(x)p_0(x)\rmd x.
\end{equation}
The equations for \(u\) and \(p\) imply
\begin{align*}
    \frac{\rmd}{\rmd t}\int_\Omega u_tp_t\rmd x
    &=
    \int_\Omega p_t\,\partial_tu_t\rmd x
    -\int_\Omega
        \nabla u_t\cdot\bigl(\nabla f(p_t)+p_t \nabla (V+W*p_t))\bigr)\rmd x\\
    &=
    -\int_\Omega
        \bigl(f(p_t)\Delta u_t+\nabla f(p_t)\cdot\nabla u_t\bigr)\rmd x\\
    &=
    -\int_{\partial\Omega}f(p_t)\nabla u_t\cdot n\rmd S
    =0.
\end{align*}
Integrating in time and using \eqref{eq:expectation}, we obtain
\[
    \int_\Omega\varphi(x)p_T(x)\rmd x
    =\int_\Omega u_0(x)p_0(x)\rmd x
    =\mathbb{E}\bigl[\varphi(X_T)\bigr],
    ~ \varphi\in C_c^\infty(\Omega).
\]
%Thus the restriction of \(\operatorname{Law}(X_T)\) to \(\Omega\) is \(p_T(x)\,dx\).
Since
\[
    \operatorname{Law}(X_T)(\Omega)
    =
    \int_\Omega p_T(x)\rmd x
    =1
\]
and \(X_T\in\overline\Omega\) almost surely, it follows that
\[
    \mathbb{P}(X_T\in\partial\Omega)=0.
\]
Hence, as probability measures on \(\overline\Omega\),
\[
    \operatorname{Law}(X_T)=p_T(x)\rmd x.
\]
The arbitrary of \(T>0\) proves
\[
    \operatorname{Law}(X_t)=p_t(x)\rmd x,
    ~ t\ge0.
\]
\end{proof}

\subsection{Proof of Theorem \ref{thm:diss} and Corollary \ref{cor:diss}}\label{proof-thm-diss} 
\begin{proof}[Proof of Theorem \ref{thm:diss}]
Applying It\^o's formula to \eqref{eq:reflected-SDE}, we obtain
\[
\begin{aligned}
\rmd \theta_t(X_t)
={}&\partial_t \theta_t(X_t)\rmd t
-\big\langle\nabla \theta_t(X_t),\nabla V(X_t)+\nabla W*p_t(X_t)\big\rangle\rmd t\\
&+\frac{f(p_t(X_t))}{p_t(X_t)}\Delta \theta_t(X_t)\rmd t
+\sqrt{\frac{2f(p_t(X_t))}{p_t(X_t)}}
\big\langle\nabla \theta_t(X_t),\rmd B_t\big\rangle\\
&-\big\langle\nabla \theta_t(X_t),n(X_t)\big\rangle\rmd L_t.
\end{aligned}
\]
The definition \eqref{def-D} gives
\begin{equation}\label{114-1}
\begin{aligned}
\rmd\theta_t(X_t)
={}&D_t(X_t)\rmd t
+\sqrt{\frac{2f(p_t(X_t))}{p_t(X_t)}}\big\langle\nabla \theta_t(X_t),\rmd B_t\big\rangle\\
&-\partial_n \theta_t(X_t)\rmd L_t.
\end{aligned}
\end{equation}
Using \eqref{def-v} and \eqref{h-Phi}, we have
\begin{equation}\label{proof:grad-v}
\nabla \theta_t(x)=\frac{f'(p_t(x))}{p_t(x)}\nabla p_t(x)
+\nabla V(x)+\nabla W*p_t(x).
\end{equation}
For $x\in\partial\Omega$,
\[
\partial_n \theta_t(x)=\frac{1}{p_t(x)}
\big[\nabla f(p_t(x))+p_t(x)\big(\nabla V(x)+\nabla W*p_t(x)\big)\big]\cdot n(x).
\]
The boundary condition in \eqref{eq:1.1} implies
\begin{align}\label{proof:normal-v}
\partial_n \theta_t(x)=0,~ (t,x)\in[0,T]\times\partial\Omega.
\end{align}
Since $L_t$ increases only when $X_t\in\partial\Omega$ by
\eqref{eq:reflected-SDE}, we get
\[
\int_0^t\partial_n \theta_s(X_s)\rmd L_s=0.
\]
Integrating \eqref{114-1} over $[0,t]$ and using the definitions of $F_t$, $M_t$ in \eqref{Mt-Ft}, we
arrive at \eqref{v-decomposition}.

By \eqref{Mt-Ft}, It\^o's isometry, Lemma \ref{PDEsolution} and Lemma \ref{SDEsolution}, we have
\begin{equation}\label{M2-proof}
\begin{aligned}
\sup_{[0,T]}\E[|M_t|^2]
&=2\sup_{[0,T]}\E\Big[\int_0^t
\frac{f(p_s(X_s))}{p_s(X_s)}
|\nabla \theta_s(X_s)|^2\rmd s \Big]\\
&=2\sup_{[0,T]}\int_0^t\int_\Omega
f(p_s(x))|\nabla \theta_s(x)|^2\rmd x\rmd s\\
&\leq 2\int_0^T\int_\Omega
f(p_s(x))|\nabla \theta_s(x)|^2\rmd x\rmd s
<+\infty.
\end{aligned}
\end{equation}
Thus $M_t$ is an $L^2$-bounded martingale.

By Lemma \ref{SDEsolution} and \eqref{def-D},
\begin{equation}\label{proof:mean-D}
\begin{aligned}
\E[D_t(X_t)]
={}&\int_\Omega p_t(x)\partial_t \theta_t(x)\rmd x
+\int_\Omega f(p_t(x))\Delta \theta_t(x)\rmd x\\
&-\int_\Omega p_t(x)\big\langle\nabla \theta_t(x),
\nabla V(x)+\nabla W*p_t(x)\big\rangle\rmd x.
\end{aligned}
\end{equation}
Differentiating \eqref{def-v}, using \eqref{h-Phi} and the symmetry
of $W$ in Assumption \ref{potential}, we get
\[
\begin{aligned}
\partial_t \theta_t(x)
={}&\frac{f'(p_t(x))}{p_t(x)}\partial_t p_t(x)
+\int_\Omega W(x-y)\partial_t p_t(y)\rmd y\\
&-\int_\Omega f'(p_t(y))\partial_t p_t(y)\rmd y
-\int_\Omega\int_\Omega W(y-z)\partial_t p_t(y)p_t(z)\rmd z\rmd y.
\end{aligned}
\]
Multiplying both sides by by $p_t(x)$, integrating and
using \eqref{inte-p}, we obtain
\begin{equation}\label{proof:time-zero}
\begin{aligned}
& \int_\Omega p_t(x)\partial_t \theta_t(x)\rmd x \\
=&\int_\Omega f'(p_t(x))\partial_t p_t(x)\rmd x
+\int_\Omega\int_\Omega W(x-y)\partial_t p
_t(y)p_t(x)\rmd y\rmd x\\
&-\int_\Omega f'(p_t(y))\partial_t p_t(y)\rmd y
-\int_\Omega\int_\Omega W(y-z)\partial_t p_t(y)p_t(z)\rmd z\rmd y=0,
\end{aligned}
\end{equation}
where the double integrals cancel by Fubini's theorem and
Assumption \ref{potential}.
Integration by parts and \eqref{proof:normal-v} give
\begin{equation}\label{proof:parts}
\begin{aligned}
&\int_\Omega f(p_t(x))\Delta \theta_t(x)\rmd x \\
=&\int_{\partial\Omega}f(p_t(x))\partial_n \theta_t(x)\rmd S(x)
-\int_\Omega\big\langle\nabla f(p_t(x)),\nabla \theta_t(x)\big\rangle\rmd x\\
=&-\int_\Omega\big\langle\nabla f(p_t(x)),\nabla \theta_t(x)\big\rangle\rmd x.
\end{aligned}
\end{equation}
Substituting \eqref{proof:time-zero} and \eqref{proof:parts}
into \eqref{proof:mean-D}, we have
\begin{equation}\label{114-3}
\begin{aligned}
& \E[D_t(X_t)] \\
=&-\int_\Omega\big\langle\nabla f(p_t(x)),\nabla \theta_t(x)\big\rangle\rmd x \\
&-\int_\Omega p_t(x)\big\langle\nabla \theta_t(x),
\nabla V(x)+\nabla W*p_t(x)\big\rangle\rmd x\\
=&-\int_\Omega\Big\langle\nabla f(p_t(x))
+p_t(x)\bigl(\nabla V(x)+\nabla W*p_t(x)\bigr),
\nabla \theta_t(x)\Big\rangle\rmd x.
\end{aligned}
\end{equation}
By \eqref{proof:grad-v} and the chain rule,
\[
\begin{aligned}
p_t(x)\nabla \theta_t(x)
={}&f'(p_t(x))\nabla p_t(x)
+p_t(x)\bigl(\nabla V(x)+\nabla W*p_t(x)\bigr)\\
={}&\nabla f(p_t(x))
+p_t(x)\bigl(\nabla V(x)+\nabla W*p_t(x)\bigr).
\end{aligned}
\]
Then by \eqref{114-3}, \eqref{h-Phi} and \eqref{def-I}, we obtain
\begin{align}\label{proof:mean-D-I}
\E[D_t(X_t)]
&=-\int_\Omega\big\langle p_t(x)\nabla \theta_t(x),\nabla \theta_t(x)\big\rangle\rmd x\notag\\
&=-\int_\Omega p_t(x)|\nabla \theta_t(x)|^2\rmd x\\
&=-\int_\Omega p_t(x)
\bigl|\nabla \psi(p_t(x))+\nabla V(x)+\nabla W*p_t (x)\bigr|^2\rmd x\notag\\
&=-\I(p_t).\notag
\end{align}
Applying Fubini's theorem, \eqref{Mt-Ft} and \eqref{proof:mean-D-I}, we finally conclude that 
\[
\E[F_t]=-\int_0^t\I(p_s)\rmd s>-\infty,
\]
where finiteness follows from Lemma \ref{PDEsolution},
Assumption \ref{potential} and \eqref{def-I}.
This proves \eqref{EFt}.
\end{proof}

\begin{proof}[Proof of Corollary \ref{cor:diss}]
By Lemma \ref{SDEsolution}, \eqref{inte-p} and \eqref{def-v}, we have
\[
\begin{aligned}
\E[\theta_t(X_t)]
={}&\int_\Omega p_t(x)\psi(p_t(x))\rmd x
+\int_\Omega V(x)p_t(x)\rmd x\\
&+\int_\Omega W*p_t(x)p_t(x)\rmd x
-\int_\Omega f(p_t(y))\rmd y\\
&-\frac12\int_\Omega W*p_t(y)p_t(y)\rmd y.
\end{aligned}
\]
Using \eqref{Phi-identity} and \eqref{def-F}, we obtain
\begin{equation}\label{proof:Ev-F}
\begin{aligned}
\E[\theta_t(X_t)]
=&\int_\Omega\Psi(p_t(x))\rmd x
+\int_\Omega V(x)p_t(x)\rmd x\\
&+\frac12\int_\Omega W*p_t(x)p_t(x)\rmd x\\
=&\F(p_t).
\end{aligned}
\end{equation}
By Theorem \ref{thm:diss}, $M_t$ is a martingale, which means
\[\E[M_t-M_{t_0}]=0.\]
Then subtracting \eqref{v-decomposition} at $t_0$ from
\eqref{v-decomposition} at $t$ and taking expectations, we get
\[
\E[\theta_t(X_t)]-\E[\theta_{t_0}(X_{t_0})]
=\E[F_t]-\E[F_{t_0}].
\]
Applying \eqref{proof:Ev-F} and \eqref{EFt},
we obtain \eqref{energy-dissipation}.
Lemma \ref{PDEsolution} and \eqref{def-I} imply
\[s\mapsto\I(p_s) \text{ is continuous on } [0,T].\]
Differentiating \eqref{energy-dissipation} with respect to $t$
at $t=t_0$ proves \eqref{eq:diss-rate}.

By Theorem \ref{thm:diss},
\[\E[M_t-M_{t_0}|\F_{t_0}]=0.\]
Taking conditional expectations of the increment in
\eqref{v-decomposition} and using \eqref{Mt-Ft}, we obtain that
\[
\E[\theta_t(X_t)\mid\F_{t_0}]-\theta_{t_0}(X_{t_0})
=\E\left[
\int_{t_0}^t D_s(X_s)\rmd s
\,\middle|\,\F_{t_0}
\right].
\]
Since $D_{t_0}(X_{t_0})$ is $\F_{t_0}$-measurable,
the $L^1$ contraction property of conditional expectation gives
\begin{equation}\label{proof:conditional-estimate}
\begin{aligned}
&\E\left[
\left|
\frac{\E[\theta_t(X_t)\mid\F_{t_0}]-\theta_{t_0}(X_{t_0})}{t-t_0}
-D_{t_0}(X_{t_0})
\right|
\right]\\
=&
\E\left[
\left|
\E\left[
\frac{1}{t-t_0}\int_{t_0}^t
\bigl(D_s(X_s)-D_{t_0}(X_{t_0})\bigr)\rmd s
\,\middle|\,\F_{t_0}
\right]
\right|
\right]\\
\leq &
\frac{1}{t-t_0}\int_{t_0}^t
\E\bigl[|D_s(X_s)-D_{t_0}(X_{t_0})|\bigr]\rmd s.
\end{aligned}
\end{equation}
By Lemmas \ref{PDEsolution}, \ref{SDEsolution} and
\eqref{def-D}, the process $D_s(X_s)$ is bounded and
has continuous paths. Dominated convergence theorem implies
\[
\lim_{s\downarrow t_0}
\E\bigl[|D_s(X_s)-D_{t_0}(X_{t_0})|\bigr]=0.
\]
Therefore,
\begin{equation}\label{proof:D-average-limit}
\lim_{t\downarrow t_0}
\frac{1}{t-t_0}\int_{t_0}^t
\E\big[\big|D_s(X_s)-D_{t_0}(X_{t_0})\big|\big]\rmd s=0.
\end{equation}
Combining \eqref{proof:conditional-estimate} and
\eqref{proof:D-average-limit} proves
\eqref{eq:conditional-diss} in $L^1(\mathbb P)$.
\end{proof}

\subsection{Proof of Lemma \ref{thm:perturbed-diss}, Corollary \ref{cor:perturbed-diss}, Lemma \ref{dervar} and Theorem \ref{main}}\label{W-proof}
\begin{proof}[Proof of Lemma \ref{thm:perturbed-diss}]
In this proof, \eqref{proof:grad-v}, \eqref{M2-proof},
\eqref{proof:time-zero} and \eqref{proof:parts}
are used with $(p,\theta)$ replaced by $(p^\var,\theta^\var)$
and $[0,T]$ replaced by $[t_0,T_\var]$.
By the boundary condition in \eqref{eq:1.2} and
Assumption \ref{ass-perturbation},
\[
\partial_n \theta^\var_t(x)=-\partial_n\var(x)=0,
~ x\in\partial\Omega.
\]
Since $L^\var$ increases only when $X^\var\in\partial\Omega$, applying It\^o's formula to the perturbed SDE \eqref{eq:reflected-SDE-per}, we get
\[
-\int_{t_0}^t
\partial_n \theta^\var_s(X_s^\var)\rmd L_s^\var=0.
\]
Integrating the It\^o formula for $\theta^\var_t(X_t^\var)$
over $[t_0,t]$, using 
\eqref{def-D-epsilon} and \eqref{Mt-Ft-epsilon},
we obtain \eqref{v-epsilon-decomposition}.
By \eqref{M2-proof}, we have
\[
\sup_{[t_0,T_\var]}\E[|M_t^\var|^2]<+\infty.
\]
Thus $M_t^\var$ is an $L^2$-bounded martingale.

By \eqref{proof:time-zero}, \eqref{proof:parts},
and \eqref{proof:grad-v} in \eqref{def-D-epsilon},
we conclude that
\[
\begin{aligned}
\E[D^\var_t(X_t^\var)]
=&-\int_\Omega p_t^\var(x)|\nabla \theta^\var_t(x)|^2\rmd x
-\int_\Omega p_t^\var(x)
\langle\nabla \theta^\var_t(x),\nabla\var(x)\rangle\rmd x\\
=&-\I(p_t^\var)
-\E\bigl[
\langle\nabla \theta^\var_t(X_t^\var),
\nabla\var(X_t^\var)\rangle
\bigr],
\end{aligned}
\]
where the last equality uses \eqref{def-I}.
Applying Fubini's theorem to the definition \eqref{Mt-Ft-epsilon}
proves \eqref{EFt-epsilon}.
The finiteness follows from the regularity and positive
lower bound of $p_t^\var$ and the smoothness of $\var$.
\end{proof}

\begin{proof}[Proof of Corollary \ref{cor:perturbed-diss}]
The identity \eqref{proof:Ev-F}, with $(p,\theta,X)$ replaced
by $(p^\var,\theta^\var,X^\var)$, implies
\[
\E[\theta^\var_t(X_t^\var)]=\F(p_t^\var).
\]
Taking expectations in \eqref{v-epsilon-decomposition}
and using \eqref{EFt-epsilon} proves
\eqref{eq:perturbed-energy-dissipation}.

Differentiating \eqref{eq:perturbed-energy-dissipation}
at $t=t_0$ from the right, using
$p_{t_0}^\var=p_{t_0}$,
$X_{t_0}^\var=X_{t_0}$ and
\eqref{proof:grad-v}, we arrive at \eqref{eq:perturbed-diss-rate}.

The conditional expectation argument in the proof of
Corollary \ref{cor:diss} applies with $(\theta,D,X,M,F)$
replaced by $(\theta^\var,D^\var,X^\var,M^\var,F^\var)$.
In particular, by 
\eqref{proof:conditional-estimate} and 
\eqref{proof:D-average-limit}, we finally obtain
\eqref{eq:conditional-perturbed-diss} in $L^1(\mathbb P)$.
\end{proof}

\begin{proof}[Proof of Theorem \ref{main}]
By the definitions of \eqref{eq:1.1} and \eqref{eq:1.2}, the perturbed and unperturbed systems coincide at any fixed time $t_0$. In particular, we have $p_{t_0}=p^{\var}_{t_0}$ and $X_{t_0}=X^{\var}_{t_0}$. 
Substituting \eqref{eq:diss-rate} and \eqref{mder} into the expression for the Wasserstein metric slope, we obtain
\begin{equation*}
\big|\partial \mathcal{F} \big|_{W_2}(p_{t_0}) = 
-\big\| 
\nabla v_{t_0}(X_{t_0})
\big\|_{L^2(\P)} \leq 0.
\end{equation*}
Similarly, using Corollary \ref{cor:perturbed-diss} and  Lemma \ref{dervar}, we obtain
\begin{equation*}
\big|\partial \mathcal{F} \big|_{W_h}(p^{\var}_{t_0}) = -\bigg \langle \nabla v^\var_{t_0}(X_{t_0}),  \dfrac{\nabla v^\var_{t_0}(X_{t_0})+\nabla \var(X_{t_0}) }{\big\| \nabla v^\var_{t_0}(X_{t_0})+\nabla\var(X_{t_0}) \big\|_{L^2(\P)}}
\bigg \rangle_{L^2(\P)}.
\end{equation*}
The Cauchy-Schwarz inequality implies
\[
\big|\partial \mathcal{F} \big|_{W_2}(p^{\var}_{t_0})  \geq \big|\partial \mathcal{F} \big|_{W_2}(p_{t_0}).
\]
Equality holds if and only if $v_t=\kappa\nabla \var$ for some constant $\kappa>0$. 

By the arbitrariness of $t_0$ and $\var$, equation \eqref{gradient} can be rewritten in the gradient flow form, which is equivalent to \eqref{form2}.
\end{proof}

\section*{Acknowledgements}
% We are grateful to Professor Nicola Gigli for insightful comments.
This work is supported by National Key R\&D Program of China (No. 2023YFA1009200), NSFC (Grants 12531009 and 11925102).

\end{document}